\input diagrams.tex  

\def\ii{i}
\def\jj{j}

\centerline{\tenbf Homotopy equivalence of isospectral graphs.}

\vskip .25 true cm
\centerline{Terrence Bisson \quad\& \quad Aristide Tsemo}
\centerline{bisson@canisius.edu \quad \  \quad tsemo58@yahoo.ca}
\vskip .25 true cm

\centerline{Abstract.}
 \bigskip
In this paper, we investigate the Quillen  model structure defined by Bisson and Tsemo in  the category of directed graphs Gph.  In particular, we give a precise description of the homotopy category  of graphs 
associated to this model structure.  We endow the categories  of N-sets and  Z-sets  with  related model structures, and show that their homotopy categories are  Quillen equivalent to the homotopy category Ho(Gph). This enables us to show that Ho(Gph) is equivalent to the category cZSet of periodic Z-sets, and to show that two finite directed graphs are almost-isospectral if and only if they are homotopy-equivalent in our sense.

\beginsection ${\cal x}$0. Introduction.

Mathematicians often study complicated categories by means of {\it invariants} 
(which are equal for isomorphic objects in the category).
Sometimes a complicated category can be replaced by 
a (perhaps simpler) homotopy category which is 
better related to the various invariants used to study it.
In topology, this was first achieved by declaring two continuous functions 
to be equivalent when one could be deformed into the other.  
But it was eventually realized that most of the important features of this analysis
are determined by the class of homotopy equivalences in the category.

Quillen [1967] presented an abstraction of this method
that applies to many categories.  A Quillen model structure on a category ${\cal E}$
works with three classes of morphisms in the category, which are assumed to
satisfy certain axioms.
Quillen described the associated homotopy category 
${\rm Ho}({\cal E})$ as a {\it localization} or 
{\it category of fractions} with respect to the class of {\it weak equivalences}
for the model structure; he defined the morphism sets for this
homotopy category by using the classes of {\it fibrations} and
{\it cofibrations} for the model structure.  We review this in section 2 here.

In Bisson and Tsemo [2008] we gave a model structure for a particular
category Gph of directed and possibly infinite graphs, 
with loops and multiple arcs allowed
(we give a precise definition of Gph in section 1 here).
We focussed on invariants in Gph defined in terms of cycles, and
defined the weak equivalences for our model structure to be the
Acyclics (graph morphisms which preserve cycles).
The cofibrations and fibrations for the model
are determined from the class of Whiskerings 
(graph morphisms produced by grafting trees).
We review this model structure in section 2 here.

The main goal of the present paper is to prove that the homotopy
category Ho(Gph) for our model structure is equivalent to 
the category cZSet of periodic Z-sets.  The proof is in section 4 here.  This
result is applied in section 5 to show that isospectral and almost-isospectral finite graphs are homotopy equivalent for our model structure.

We use  the fact that  whiskered cycles (disjoint unions of cycles with
trees attached to them) are cofibrant objects in our model structure.
These graphs can also be described as Cayley graphs of N-sets, where N is the monoid of the 
natural numbers under addition (with 1 as generator).  We make even more use of 
disjoint union of cycles as Cayley graphs of Z-sets.

Each of the categories Gph, NSet, and ZSet is a presheaf topos, 
with adjoint functors relating them.   By selecting
appropriate adjoint functors, we transport our model structure from Gph
to the categories NSet and ZSet.  Most of the functors used here and throughout the 
paper arise in groups of three $(F,G,H)$, made up of two overlapping adjunctions $(F,G)$ and
$(G,H)$ between
presheaf categories.  See the end of section 1 for background.

We show that Ho(ZSet) and Ho(NSet) are both equivalent to Ho(Gph).
We do this by using a further adjunction between ZSet and the category
cZSet of periodic Z-sets.  In fact, we exhibit Quillen equivalences between 
cZSet and ZSet, NSet, and Gph, where we use a trivial model 
structure on cZSet.

\smallskip
Here is a more detailed outline of the sections of this paper.

In section 1 we define the category Gph, which is a presheaf category, and thus a topos.
We define a subcategory NGph of Gph, which is equivalent
to the category NSet of actions of the additive monoid of natural numbers. 
We note that NGph is also a topos.
Then we give a similar discussion of ZGph (which is equivalent to the topos ZSet
of actions of the additive group of integers), and TGph (which is 
equivalent to the topos Set).
We observe, in passing, that these equivalences provide (very simple) examples
of Grothendieck's version of Galois theory.  
We show that these subcategories are reflective and coreflective subcategories of Gph.   
The functors we need for this  are arising as
``adjoint triples'' $(F,G,H)$ of functors between presheaf categories. We establish our
our conventions for these functors at the end of section 1; they are used almost
everywhere in the paper.

In section 2 we recall the definitions of model structure and
cofibrantly-generated model structure.  We recall our terminology
of Surjecting, Whiskering, and Acyclic graph morphisms from Bisson and Tsemo [2008],
and the definition of our model structure on Gph.  
We show that our model structure is cofibrantly-generated, and
 use a general (folk) theorem on the
transport of cofibrantly-generated model structures to define
model structures on NSet (and NGph), and ZSet (and ZGph).

In section 3, we analyze these new model structures on NSet and ZSet,
with especial attention to fibrant objects and cofibrant objects.
Motivated by this analysis, we develop a cofibrant replacement functor for
the category Gph.  Our construction uses the coreflection functor $H$ for ZGph
as a subcategory of Gph.

In section 4 we give some background on homotopy functors and on
Quillen's construction of the homotopy category ${\rm Ho}({\cal E})$ associated to 
a model structure on ${\cal E}$, and on his construction of
derived functors (for adjoint functors satisfying appropriate conditions).  
We look for examples of homotopy functors on Gph,
and for examples which satisfy the Quillen adjunction conditions.
We use a particular adjunction relating Gph and ZSet
(with left adjoint the functor $H$ which assigns to each graph the set of
all its bi-infinite paths)  to show that
Ho(Gph) is equivalent to the category cZSet of periodic Z-sets.

In section 5 we use the functor $H$ to
associate a zeta series $Z_X(u)$
to each almost-finite graph $X$.  This fits very well with work
of Dress and Siebeneicher [1988] on 
the Burnside ring of the category  of almost-finite Z-sets.  
As a consequence of our calculation of Ho(Gph), we show that finite graphs are 
almost-isospectral if and only if they
are homotopy equivalent (that is, isomorphic in the homotopy category).

\beginsection ${\cal x}$1. Some subcategories of graphs. 

This paper is about Gph, a convenient category of
graphs, precisely described in the paragraph below. 
In Bisson and Tsemo [2008] we introduced a Quillen model structure on Gph.
Here we will show how to study that structure, 
and the resulting homotopy category, 
by means of some of its subcategories.

We define a {\it graph} to be a data-structure $X=(X_0,X_1,s,t)$ with a set $X_0$ of {\it nodes},
a set $X_1$ of {\it arcs}, and a pair of functions $s,t:X_1\to X_0$ which
specify the {\it source} and {\it target} nodes of each arc.
We may say that $a\in X_1$ is an arc which {\it leaves node $s(a)$ and enters node $t(a)$};
and that a {\it loop} is an arc $a$ with $s(a)=t(a)$.
A {\it graph morphism} $f:X\to Y$ is a pair of functions 
$f_1:X_1\to Y_1$ and $f_0:X_0\to Y_0$ such that
$s\circ f_1=f_0\circ s$ and $t\circ f_1=f_0\circ t$.
This defines the particular category Gph that we study here.

In fact,  Gph is the category of presheafs on a small category; see
Lawvere [1989]  or Lawvere and Schanuel [1997] for fascinating discussions.
It follows that ${\rm Gph}$ is a topos, and thus a category with many nice geometric and
algebraic and logical properties; see Mac Lane and Moerdijk [1994], for instance.

In this paper we want to consider some very special kinds of graphs, as follows.

\medskip\noindent{\bf Definition:} A graph $X$ is

1): an {\it N-graph} when each node of $X$ has exactly one arc entering.  

2): a {\it Z-graph} when each node of $X$ has exactly one arc entering and   
exactly one arc leaving.

3): a {\it T-graph} when each node of $X$ has exactly one loop, 
and $X$ has no other arcs.    

A T-graph might be called a {\it terminal} graph (or {\it graph of loops})
and a Z-graph might be called a {\it graph of cycles}.  
An N-graph might be called a {\it graph of whiskered cycles}.

Let NGph denote the {\it full subcategory} of Gph whose objects are the N-graphs;
this means that we take all graph morphisms between N-graphs
as the morphisms in NGph.
Similarly, let ZGph denote the full subcategory of Gph whose objects are the Z-graphs,
and let TGph denote the full subcategory of Gph whose objects are the terminal graphs.
We have a chain of subcategories
$$\hbox{TGph}\subset \hbox{ZGph}\subset \hbox{NGph} \subset {\rm Gph}$$
In fact, these Gph subcategories are equivalent to some well-known categories.
Let $G$ be a monoid, with associative binary operation $G\times G\to G:(g,h)\mapsto g*h$ 
and with neutral element $e$; 
a $G$-set is a set $S$ together with an action $\mu:G\times S\to S$ 
such that $\mu(e,x)=x$ and $\mu(g,\mu(h,x))=\mu(g*h,x)$.

Consider the monoid $N$ of natural numbers under addition,
and the group $Z$ of integers under addition.  
A set $S$ together with an arbitrary function $\sigma:S\to S$
defines an action by $\mu(n,x)=\sigma^n(x)$ for $n\in N$.
A set $S$ together with an arbitrary {\it invertible} function $\sigma:S\to S$
defines an action by $\mu(n,x)=\sigma^n(x)$ for $n\in Z$.
So in this paper we will use the following alternative
definitions of (the categories of) N-sets and Z-sets.

\medskip\noindent{\bf Definition:} 
Let Set denote the category of sets.
Let NSet denote the category of N-sets; here
an N-set is a pair $(S,\sigma)$ with $\sigma$ a function from $S$ to $S$; and
a map of N-sets from $(S,\sigma)$ to $(S',\sigma')$ is a function $f:S\to S'$ such that
$\sigma'\circ f=f\circ\sigma$.  
Let ZSet denote the full subcategory of NSet with objects $(S,\sigma)$ 
where $\sigma$ is a bijection.

\medskip
For any N-set $(S,\sigma)$, we define a graph $X=G(S,\sigma)$ with nodes $X_0=S$ and arcs $X_1=S$, 
and with $s,t:X_1\to X_0$ given by $s(x)=\sigma(x)$  and $t(x)=x$ for each $x\in S$.  
Thus the elements in the N-set $S$ give the nodes and the arcs in the graph $X$,
and each arc $x$ has source $\sigma(x)$ and target $x$.
In the N-set $S$ we think of $\sigma(x)$
as telling the unique ``source'' or ``parent'' of each element $x$.

Note that we are directing our arcs opposite to the way that seems natural
in graphical representation of dynamical systems (see Lawvere and Schanuel [1997], for instance).
But our convention is designed to fit well with
the notion of ``whiskerings'' in our model structure
(see section 2 here).
  
\medskip\noindent{\bf Proposition:} $G$ is a functor from NSet to Gph; moreover,

a) the functor $G$ gives an equivalence from the category NSet 
to the sub-category NGph in Gph; and

b) the restriction of $G$ gives an equivalence from the category
ZSet to the sub-category ZGph in Gph.

c) the restriction of $G$ gives an equivalence from the category
Set to the sub-category TGph in Gph.

\smallskip\noindent{\bf Proof:}  If $f:(S,\sigma)\to (S',\sigma')$ is a map of N-sets, 
we define a graph morphism 
$G(f):G(S,\sigma)\to G(S',\sigma')$
by $G(f)_0(x)=G(f)_1(x)=f(x)$ for $x$ in $S$.  This preserves composition.
We note that $G(S,\sigma)$ is an N-graph, and every
N-graph has a unique isomorphism to a graph $X$ in the image of $G$ 
(where $X_0=X_1$ and $t$ is the identity). 
If $X$ is in NGph we define an N-set
$H(X)=(X_0,\sigma)$ by $\sigma(x)=s(a)$ where $a$ is the unique arc entering the node $x$,
and a graph morphism $g:X\to Y$ gives $H(g):H(X)\to H(Y)$ by $H(g)(x)=g_0(x)$ for $x\in X_0$.
This preserves composition and gives a functor from NGph to NSet.
Thus we have 
$$G:\hbox{NSet}\to \hbox{NGph}\quad{\rm and}\quad H:\hbox{NGph}\to \hbox{NSet}.$$
Note that $H(G(S,\sigma))=(S,\sigma)$ and $G(H(X))=X$.
In fact, $G$ and $H$ give inverse bijections between the set of N-set maps 
$(S,\sigma)$ to $(S',\sigma')$ and the 
set of graph morphisms from $G(S,\sigma)$ to $G(S',\sigma')$.
This says that the functor $G$ is full and faithful, with image the category NGph.
We can carry out a similar analysis for the restriction of $G$ to the sub-category ZSet, and
the restriction of $H$ to the sub-category ZGph. 
The analysis for Set and TGph is also similar (and rather trivial).
QED.

\smallskip
Each part of the above proof exhibits an ``adjoint pair'' of functors $(G,H)$
(as discussed below), and shows that it gives an equivalences of categories.
The proof can also be understood as an example of (the representable case of)
Grothendieck's Galois theory.  The paper by Dubuc and de la Vega [2000]
gives a self-contained exposition which seems relevant to our examples here.

It is interesting to note that the calculation of products is the same at each level of
the inclusions 
$$\hbox{TGph}\subset\hbox{ZGph}\subset\hbox{NGph}\subset\hbox{Gph}$$
The same is true of coproducts.  
In fact, the functor $G:\hbox{NSet}\to\hbox{Gph}$ preserves all limits and colimits,
since it has left and right adjoints; and similarly for the other $G$ functors.
These left and right adjoints show that TGph and NGph and ZGph are 
``(full) reflective and coreflective subcategories'' of Gph.  
We will use this in Section 2 to  transport a Quillen model structure from Gph
to NGph and ZGph.  Let us explain this terminology, starting
with the case TGph of terminal graphs, which is especially simple.

Recall that  TGph denotes the full subcategory of terminal graphs,
those graphs which are disjoint unions of $1$.
For each graph $X$ we will describe a graph morphism
$\sqcap_X:X\to T_X$ which is ``universal'' among morphisms from $X$ to terminal graphs, 
in that any $f:X\to T'$ with $T'$ terminal factors 
through a unique graph morphism $f':T_X\to T'$. 
Dually, we will also describe a graph morphism
$\sqcup_X:T''_X\to X$ which is ``universal'' among morphisms from terminal graphs to $X$.
The existence of these natural, universal examples
means (loosely speaking) that TGph is a ``reflective and coreflective subcategory'' of Gph.
These are best described in the language of adjoint functors.  Let us give a
quick review of some standard definitions (see Mac Lane [1971], for instance).

An {\it adjunction} between categories ${\cal X}$ and ${\cal Y}$
is a pair $(L,R)$ of functors $L:{\cal X}\to {\cal Y}$ and $R:{\cal Y}\to {\cal X}$ 
together with a natural bijection of morphism sets 
${\cal Y}(L(X),Y)\to{\cal X}(X,R(Y))$.
In this case, we may say that $(L,R)$ is an {\it adjoint pair}, with $L$ as the {\it left adjoint}
and $R$ as the {\it right adjoint}, and denote this by
$$L:{\cal X}\rightleftharpoons {\cal Y}:R$$
For any $(L,R)$ is an adjoint pair we have a natural transformation $X\to RL(X)$, called the
{\it unit} of the adjunction; dually, there is a
natural transformation $LR(Y)\to Y$, called the
{\it counit} of the adjunction.  In the rest of this section, we use two
special types of adjunctions.

\smallskip\noindent
{\bf Definition:}
A subcategory ${\cal E}'$ of ${\cal E}$  is a {\it reflective subcategory} 
when the inclusion functor $G: {\cal E}'\to {\cal E}$ has a left adjoint functor $F$,
with adjunction $(F,G)$; then $F$ is the {\it reflection} functor. 
Dually,  ${\cal E}'$ is a {\it coreflective subcategory} of ${\cal E}$ when
the inclusion functor $G: {\cal E}'\to {\cal E}$ has a right adjoint functor $H$,
with adjunction $(G,H)$; then  $H$ is the {\it coreflection} functor.

\smallskip
As an example, let us show that TGph is reflective and coreflective in Gph.
In order to show that TGph is a reflective subcategory of Gph, we use the adjunction
$$F:\hbox{Gph}\rightleftharpoons\hbox{Set}:G$$
Here $G$ is the functor from Set to Gph which assigns to set $S$ the terminal graph 
with one loop for each element of $S$,
and $F(X)$ is the set of components of the graph $X$.
This can be defined as the set
of equivalence classes of nodes of $X$, with respect to the equivalence relation generated by 
the source and target functions $s,t:X_1\to X_0$.
The unit of the adjunction $X\to GF(X)$ is universal
among graph morphisms from $X$ to terminal graphs,
as mentioned above.
This shows the desired the adjunction.
We may use the notation  $F(X)=\pi_0(X)$ and $G(S)=\sum_S 1$.
The image of the functor $GF:{\rm Gph} \to {\rm Gph}$ 
is equivalent to the subcategory TGph of Gph
(in fact, any terminal graph has a unique isomorphism
to a graph $X$ with $X_0=X_1$ and $s$ and $t$ as the identity).
The counit  $v:FG(S)\to S$ of the adjunction is an isomorphism for every set $S$.


\smallskip
Dually, consider the adjunction 
$$G:\hbox{Set}\rightleftharpoons\hbox{Gph}:H$$
where $H(X)=[1,X]$ is the set of graph morphisms from $1$ to $X$, the set of those arcs of $X$
which are loops. The counit of the adjunction $GH(X)\to X$ is universal
among morphisms from terminal graphs to $X$.
This shows the adjunction. 

Let us extend the above discussion to handle the categories ZGph and NGph.

\smallskip\noindent
{\bf Proposition:} TGph and NGph and ZGph are reflective and coreflective subcategories of Gph
$$\hbox{TGph}\subset \hbox{ZGph}\subset \hbox{NGph} \subset {\rm Gph}$$

\smallskip\noindent
{\bf Proof:} We have shown that  TGph is a full reflective subcategory of Gph.
Now we give a similar treatment of NGph.  First we define the functor $F$ in the adjunction
$$F:\hbox{Gph}\rightleftharpoons \hbox{NSet}:G$$
Let ${\bf P}$ denote the unending path graph (or rooted tree); the nodes of ${\bf P}$ are
the natural numbers, and there is one arc $(n):n\to n+1$ for each $n\geq 0$.
Note that ${\bf P}$ is not itself in NGph, since the root $0$ has no arc entering it. 
We will use the graph ${\bf P}$ to define the reflection from Gph to NGph.
For any graph $X$, let $F(X)$ denote the N-set $(\pi_0({\bf P}\times X),\sigma)$.
This is an N-set since the assignment $\sigma([n,x])=[n+1,x]$ gives a well-defined  function 
on the set of connected components
of the graph ${\bf P}\times X$, since any arc $((n),a):(n,x)\to (n+1,y)$
is accompanied by an arc $((n+1),a):(n+1,x)\to (n+2,y)$.
We could instead define a graph morphism $\sigma:({\bf P} \times X)\to ({\bf P} \times X)$,
which then gives a function $\sigma: \pi_0({\bf P} \times X)\to \pi_0({\bf P} \times X)$, by
functoriality of $\pi_0$.
Let us sketch a proof that $F$ is left adjoint to $G$, and thus show that 
NGph is a reflective subcategory of Gph.
Given any graph $X$ and any N-set $S$, any morphism $f:X\to G(S)$ in Gph 
has the property that the existence of an arc $x'\to x$ in $X$ 
implies $f(x')=\sigma(f(x))$ in $S$.  We define a function
from the nodes of ${\bf P}\times X$ to $S$ by $(n,x)\mapsto \sigma^n(f(x))$.  
This is well defined on connected components of ${\bf P}\times X$ since
any arc $(n+1,x')\to (n,x)$ in ${\bf P}\times X$
implies the existence of an arc $x'\to x$ in $X$, which implies $f(x')=\sigma(f(x))$ 
and $\sigma^n(f(x'))=\sigma^{n+1}(f(x))$ in $S$.
Thus we get a morphism $g:\pi_0({\bf P}\times X)\to S$ in NSet,
since $g(\sigma[n,x]) =g([n+1,x])=\sigma^{n+1}(f(x))=\sigma(g([n,x])$.
The inverse correspondence can be established in a similar way.
The above correspondences
could also be described in terms of the unit $X\to G(F(X))$ (in Gph) or
counit $F(G(S))\to S$ (in NSet).
We merely sketch the reflection from NGph to ZGph, which also gives 
the reflection from Gph to ZGph. Informally, this functor
combs any whisker down along the cycle it comes from.  
We could instead show that ZSet is a reflective subcategory of NSet, by
describing the functor which is left adjoint to the inclusion.   QED.

For any monoid $G$, the category $G\ $Set is a presheaf category,
and thus a topos; see Mac Lane and Moerdijk [1994], for instance.
A topos has all products, and all coproducts (sums);
it also has pull-backs (fiber products) and pushouts.
In fact, since Gph, NSet, and ZSet are presheaf categories,
these categorical constructions can be performed ``elementwise''.

As simple examples, we have that
the N-set $(1,{\rm id})$ (the one point set with its identity function) 
is a {\it terminal object} in NSet;
this means that for every N-set there is a unique N-set map $(S,\sigma)\to(1,{\rm id})$.
The empty set with its identity function is an {\it initial object} in NSet;
this means that for every N-set there is a unique N-set map $(0,{\rm id})\to(S,\sigma)$.
These objects $0$ and $1$ also provide the initial and terminal objects for ZSet.

In the above discussions of NSet 
we have actually used three functors  $(F,G,H)$ between Gph and NSet.
This is an example of an ``adjoint triple'' 
between two  presheaf  categories, coming from a functor between the site categories. 
A similar remark applies to the functors between Gph and ZSet, etc.
Here is a brief sketch of the situation.

If $C$ is a small category then the topos of presheaves on $C$,
which we may denote by $C$ Set, is the category of functors from $C^{op}$ to Set.
If $\phi:C\to D$ is a functor, then we get an adjoint triple $(\phi_!,\phi^*,\phi_*)$ 
$$  \matrix{A        &\phi^*(X)&A         \cr
		  \phi_!(A)&X        &\phi_*(A) \cr} \quad{\rm with}\quad
\phi_!:C {\rm\ Set}\rightleftharpoons D {\rm\ Set}:\phi^*\quad{\rm and}\quad
\phi^*:D {\rm\ Set}\rightleftharpoons C {\rm\ Set}:\phi_*$$
This is meant to schematically display one functor
$\phi^*$ from $D$ Sets to $C$ Sets,
and two functors $\phi_!$ and $\phi_*$ from $C$ Sets to $D$ Sets.		  
The three functors here are adjoint in the sense that 
$A\to \phi^*(X)$ in $C$ Sets corresponds to  $\phi_!(A)\to X$ in $D$ Sets,
and $\phi^*(X)\to A$ in $C$ Sets corresponds to $X\to \phi_*(A)$ in $D$ Sets.
See the analysis in Expose I.5 of Grothendieck [1972]
This concept is related to that of  ``essential geometric morphism'' 
 $\phi:C {\rm\ Set}\Rightarrow D {\rm\ Set}$  in topos theory;
 see Mac Lane and Moerdijk [1994], for instance.
Almost all the adjunctions used in this paper come from such adjoint triples 
$(F,G,H)=(\phi_!,\phi^*,\phi_*)$.

\beginsection ${\cal x}$2. Quillen model structures.

Let us start with some convenient notation.

\smallskip\noindent
{\bf Definition:}
Let $\ell:X\to Y$ and $r:A\to B$ be morphisms in a category ${\cal E}$.
We say that $\ell$ is {\it weak orthogonal} to $r$ (abbreviated by $\ell\dagger r$) 
when, for all $f$ and $g$,
$${\rm if}\quad \diagram
          X & \rTo^{f} & A  \cr
   \dTo^{\ell} &           &  \dTo^{r}   \cr
          Y  & \rTo^g   &B
                                    \enddiagram \quad{\rm commutes,\ then}\quad
           \diagram
          X & \rTo^{f} & A  \cr
   \dTo^{\ell} &  \NE^{h} &  \dTo^{r}   \cr
          Y & \rTo^g   &B
                                    \enddiagram\quad {\rm commutes\ for\ some\ } h.$$ 
Given a class ${\cal F}$ of morphisms we define
${\cal F}^\dagger=\{r: f\dagger r, \ \forall f\in{\cal F}\}\quad{\rm and}\quad
{}^\dagger{\cal F}=\{\ell: \ell\dagger f, \ \forall f\in{\cal F}\}$.
A {\it weak factorization system} in ${\cal E}$ is given by two classes
${\cal L}$ and ${\cal R}$, such that ${\cal L}^\dagger={\cal R}$ and
${\cal L}={}^\dagger{\cal R}$ and such that, for any morphism $c$ in ${\cal E}$,
there exist $\ell\in{\cal L}$ and $r\in{\cal R}$ with $c=r\circ\ell$.

\medskip
Using the above, we may express Quillen's notion [1967] 
of ``model category'' via the following axioms, 
which we learned from Section 7 of Joyal and Tierney [2006].

\medskip\noindent
{\bf Definition:} Suppose that ${\cal E}$ is a category with finite limits and colimits.
A {\it model structure} on ${\cal E}$ 
is a triple $({\cal C},{\cal W},{\cal F})$
of classes of morphisms in ${\cal E}$ that satisfies

1) ``three for two'': if two of the three morphisms  $a, b, a\circ b$ belong to ${\cal W}$
then so does the third,

2) the pair $(\underline{\cal C} , {\cal F})$ is a weak factorization system 
(where $\underline{\cal C}={\cal C}\cap{\cal W}$),

3) the pair $({\cal C} , \underline{\cal F})$ is a weak factorization system 
(where $\underline{\cal F}={\cal W}\cap{\cal F}$).

\medskip
The morphisms in  ${\cal W}$ are called {\it weak equivalences}.
The morphisms in  ${\cal C}$ are called {\it cofibrations}; and
the morphisms in  $\underline{\cal C}$ are called {\it acyclic cofibrations}.
The morphisms in  ${\cal F}$ are called {\it fibrations}, and
the morphisms in  $\underline{\cal F}$ are called {\it acyclic fibrations}.

\medskip
In Bisson and Tsemo [2008] we introduced a Quillen model structure on Gph.
Its description used three types of graph morphisms, which we
called Surjectings, Whiskerings, and  Acyclics.
They can be defined as follows.

\smallskip
\item{*} A graph morphism $f:X\to Y$ is {\it Surjecting}
when the induced function $f:X(x,*)\to Y(f(x),*)$ is surjective 
for all $x\in X_0$.  Here, for any graph $Z$ and any node $z$,  $Z(z,*)$ denotes the set of arcs 
in $Z$ which have source $z$. 

\smallskip
\item{*} A graph morphism $f:X\to Y$ is {\it Acyclic} when
$C_n(f):C_n(X)\to C_n(Y)$ is bijective for all $n>0$.
Here ${\bf C}_n$ is the {\it (directed) cycle graph}, with
the integers mod $n$ as its nodes, and also as its arcs, and
with $s(i)=i+1$ and $t(i)=i$.
Then $C_n(X)$ denotes the set of graph morphisms from ${\bf C}_n$ to $X$.

\smallskip
\item{*} A graph morphism $f:X\to Y$ is a {\it Whiskering} when
$Y$ is formed by attaching rooted trees to $X$.
Here a {\it rooted tree} is a graph $T$ with a node $r$ (its {\it root})
such that, for each each node $x$ in $T$, there is a unique (directed) path in 
$T$ from $r$ to $x$.  Then ``attaching'' the rooted tree $T$ to $X$ means 
identifying the root $r$ with a node of $X$; this is forming the pushout of 
graph morphisms $r\to T$ and $r\to X$, where $r$ is considered as a graph
with one node and no arcs).

\smallskip

Here we will interpret these morphism classes, and describe our model structure for Gph,
in terms of the following standard notions (see section 2.1 in Hovey [1999], for instance).
A model structure $({\cal C},{\cal W},{\cal F})$ 
is {\it cofibrantly generated} if there are sets $I$ and $J$
of morphisms such that $J^\dagger={\cal F}$ and $I^\dagger=\underline{\cal F}$,
so that  ${\cal C}={}^\dagger(I^\dagger)$ and 
$\underline{\cal C}={}^\dagger(J^\dagger)$. 
In short, a cofibrantly-generated model structure is given by the weak factorization systems
$$({\cal C},\underline {\cal F})= ({}^\dagger(I^\dagger),I^\dagger) \quad {\rm and}\quad 
(\underline {\cal C},{\cal F})= ({}^\dagger(J^\dagger),J^\dagger).$$
We may say that $J$ generates the Acyclic Cofibrations,
and that $I$ generates the Cofibrations.  
Note that there is usually a smallness condition included in the definition,
but this smallness condition mentioned is vacuous in Gph:
every object in Gph is small with respect to every set of morphisms in Gph 
(since Gph is a presheaf category on a small category);  one can mimic the 
proof from Example 2.1.5 in Hovey, for instance.

Let us describe sets $I$ and $J$ which generate our model structure for Gph.
Let ${\bf s}:{\bf D}\to {\bf A}$ be the ``source'' graph morphism,
which exhibits the ``dot'' graph ${\bf D}$ as the source subgraph 
of the ``arrow'' graph ${\bf A}$.  More precisely, 
${\bf A}$ is the graph with two nodes, $0$ and $1$, and one arc $a$ from $0$ to $1$;
and ${\bf s}$ is the inclusion of the subgraph ${\bf D}$ with one node $0$ and no arcs.
Let ${\bf i}_n:{\bf 0}\to {\bf C}_n$ be the initial graph morphism, and let
${\bf j}_n:{\bf C}_n+{\bf C}_n\to {\bf C}_n$ be the coproduct graph morphism.
Let $J=\{{\bf s}\}$;  let $K=\{{\bf i_n}, {\bf j_n}: n>0\}$; and
let $I=J\cup K$. 

\smallskip\noindent
{\bf Theorem:}  Gph has a cofibrantly-generated model structure 
with Acyclic Cofibrations generated by $J$ and Cofibrations generated by $I$,
and with weak equivalences ${\cal W}=K^\dagger$.

\smallskip\noindent
{\bf Proof:} Here is a sketch that the morphism classes given above
satisfy the axioms for a model structure on Gph
(see Bisson and Tsemo [2008] for more details).  Consider ${\cal F}= J^\dagger$
and ${\cal C}={}^\dagger(I^\dagger)$.
The following three observations follow directly from
the definition of weak orthogonality:
1) $K^\dagger$ is precisely the Acyclic graph morphisms;
2) $J^\dagger$ is precisely the Surjecting graph morphisms; and
3) $I^\dagger=(J \cup K)^\dagger$ is precisely the Acyclic Surjecting graph morphisms.
The Surjectings form
the class ${\cal F}$ of fibrations for our model structure.
The  Acyclics form the class ${\cal W}$  of weak equivalences for our model structure.
The Acyclic Surjectings form the class 
${\cal W}\cap {\cal F} = \underline {\cal F}$ of acyclic fibrations 
for our model structure.
Then ${}^\dagger(I^\dagger)={}^\dagger({\cal W}\cap {\cal F})$.
So ${}^\dagger(I^\dagger)={\cal C}$, the cofibrations for our model structure. QED.

Note that there are general results for when ${\cal W}$ and $I$ and $J$ generate a model
structure (see section 2.1 in Hovey [1999], for instance).
The following standard notion will help us state a general theorem 
about ``transporting'' Quillen model structures to related categories.

Let ${\cal E}$ be a category with all limits and colimits. 
For a set $H$ of morphisms in ${\cal E}$, let ${\rm cell}(H)$ denote 
the class of all transfinite compositions of pushouts of elements in $H$
(see section 2.1 in Hovey [1999] for discussion of 
pushouts, transfinite compositions, retracts in the morphism category, etc).
Morphisms in ${\rm cell}(H)$ are called {\it relative $H$-cell complexes}; 
a graph $X$ is called an {\it $H$-cell complex} if $0\to X$ is a relative $H$-cell complex.  
All this is suggested by the
notion in topology of building up a space by attaching cells.

\smallskip
Here are some examples in the category Gph, where 
we take $J=\{{\bf s}\}$ and $K=\{{\bf i_n}, {\bf j_n}: n>0\}$ and $I=J\cup K$.

\smallskip\noindent
{\bf Proposition:} The morphisms in ${\rm cell}(J)$ are the Whiskerings, 
the acyclic cofibrations for our model structure.  
Moreover, ${\rm cell}(J)={}^\dagger(J^\dagger)$. 

\smallskip\noindent
{\bf Proof:} Each pushout of ${\bf s}$ attaches a single arc as a Whisker.
Attaching a rooted tree corresponds to a composition of these;
and the Whiskerings are exactly the
class of all transfinite compositions of pushouts of elements in $J$.
We also know that the Whiskerings are closed with
respect to retracts in the morphism category of Gph (see Bisson and Tsemo [2008]).
To prove the second statement, we use this, together with some general
facts from Hovey.
We always have $H\subseteq{\rm cell}(H)\subseteq{}^\dagger(H^\dagger)$.
Suppose that ``the domains of morphisms in $H$ are small with respect to ${\rm cell}(H)$''.
Then Hovey uses a general version of the small object argument 
(based on Lemma 3 of chapter II.3 in Quillen [1967] ) 
to show that  any morphism in ${}^\dagger(H^\dagger)$ is the retract,
in the category of morphisms of ${\cal E}$, of some morphism in ${\rm cell}(H)$. 
But every object in Gph is small with respect to every set of morphisms in Gph,
as mentioned above, so the  
smallness condition here is vacuous in Gph. QED.

\smallskip\noindent
{\bf Proposition:}  If $C$ is a disjoint union of  cycle graphs,
then every inclusion $X\to X+C$ is in ${\rm cell}(K)$.
Every graph morphism between disjoint unions of cycle graphs 
is in ${\rm cell}(K)$. 

\smallskip\noindent
{\bf Proof:} For any graph $X$, the morphism $X\to X+C_n$ is a pushout of $i_n$;
if $C$ is any disjoint union of cycle graphs, then $X\to X+C$ is a transfinite composition of
pushouts of $i_n$ for $n>0$.  For the second statement,
let $\pi_{n,k}:{\bf C}_{nk}\to {\bf C}_n$ (for $n>0$ and $k>0$) 
denote the graph morphism given on nodes by $\pi_{n,k}(i)=i\ {\rm mod\ }n$.
We can exhibit $\pi_{n,k}$ as a pushout of $j_{nk}$, as follows.
Consider the graph morphism $f:{\bf C}_{nk}+{\bf C}_{nk}\to {\bf C}_{nk}$ given on nodes 
by $f(i,0)=i+n$ and $f(i,1)=i$, 
where we think of graph ${\bf C}_{nk}+{\bf C}_{nk}$ as 
the product of ${\bf C}_{nk}$ and the set $\{ 0,1 \}$.
Then $\pi_{n,k}$ is the pushout of $f$ and $j_{nk}$. 
Any graph morphism between disjoint unions of cycle graphs is a pushout of such maps
(up to isomorphisms). QED.

\smallskip
Since  ${\rm cell}(K)\subseteq  {\rm cell}(I)\subseteq  {\cal C}$, these 
propositions are describing some of the cofibrations for our model structure on Gph.
But here are some morphisms which are {\bf not} cofibrations for our model structure.

\smallskip\noindent
{\bf Examples:}

a) The graph morphism $\pi_n:{\bf Z}\to {\bf C_n}$ is not a cofibration. 
We may show this by constructing an explicit  
Acyclic Surjecting graph morphism $g:X\to Y$  such that
the weak orthogonality $\pi_n \dagger g$ fails.  
Let $Y={\bf C}_n$ and let $X={\bf C}_n{\bf P}$, the graph formed by attaching
the root of the unending path ${\bf P}$  at the $0$ node in ${\bf C}_n$ (recall
that we defined ${\bf P}$ in section 1 here).
Let $f:{\bf Z}\to {\bf C}_n{\bf P}$ be the graph morphism given on nodes by
$m\mapsto m\ {\rm mod}\ n$ for $m\geq 0$ and $m\mapsto -m$ for $m\leq 0$.
Then the commutative square with horizontal arrows 
$f:{\bf Z}\to {\bf C}_n{\bf P}$ and ${\rm id}:{\bf C}_n\to {\bf C}_n$ has no lifting.
It follows that $\pi_n$ is not  in ${\rm cell}(I)$.

b) Also, $0\to {\bf Z}$ is not a cofibration.  Note that it is itself
an Acyclic Surjecting graph morphism, and it is not weakly orthogonal to itself.

c) Also, ${\bf Z}+{\bf Z}\to{\bf Z}$ is not a cofibration, since it is 
is an Acyclic Surjecting graph morphism which is not weakly orthogonal to itself.

\bigskip
We want to use our model structure on Gph to determine model structures on NSet and ZSet.
Here is a general result, referred to as
``creating model structures along a right adjoint'' (by Hirschhorn, Hopkins, Beke, etc), or
as ``transferring model structures along adjoint functors'' (by Crans, etc).
According to Berger and Moerdijk [2003] : 
``Cofibrantly generated model structures may be transferred 
along the left adjoint functor of an adjunction.
The first general statement of such a transfer in the literature is due to Crans.'' 
Here is their informal statement of this ``transfer principle''.
The reference is to Crans [1995].

\smallskip\noindent
{\bf Transport Theorem:} Let ${\cal E}$ be a model category which is cofibrantly generated,
with cofibrations generated by  $I$ and  acyclic cofibrations generated by $J$.
Let ${\cal E}'$ be a category with all limits and colimits, 
and suppose that we have an adjunction
$$L:{\cal E}\rightleftharpoons{\cal E}':R\quad{\rm with}\quad
R({\rm cell}\ L(J))\subseteq {\cal W}.$$
Also, assume that the sets $L(I)$ and $L(J)$ each permit the small object argument.
Then there is a cofibrantly generated model structure on ${\cal E}'$
with generating cofibrations $F(I)$ and generating acyclic cofibrations $F(J)$.
Moreover, the model structure $({\cal C}',{\cal W}',{\cal F}')$  satisfies
$f\in{\cal W}'$ iff $R(f)\in{\cal W}$ 
and $f\in{\cal F}'$ iff $R(f)\in{\cal F}$.

\smallskip
As mentioned before, the smallness conditions are automatically satisfied in our
presheaf categories. So, in our examples, the main hypothesis for the theorem is:
 $f\in {\rm cell}L(J)$ implies $R(f)\in{\cal W}$.

\smallskip
Let us translate some definitions from Gph into NSet.
For  an NSet map $f:(S,\sigma)\to (T,\sigma)$ we say that:

\item{a)} $f$  is {\it Acyclic} when
$C_n(f):C_n(S,\sigma)\to C_n(T,\sigma)$ is a bijection for every $n>0$.
Here $C_n(S,\sigma)=\{x\in S : \sigma^n(x)=x\}$ (we could call these the $n$-periodic points).

\item{b)} $f$ is {\it Surjecting} when
$f:\sigma^{-1}(x)\to \sigma^{-1}(f(x))$ is a surjection for every $x$ in $S$.

\item{c)} $f$ is {\it Whiskering} when
$f$ is an injective function and $x\notin f(S)$ implies that there exists some natural number
$n$ with $\sigma^n(x)\in f(S)$.

\smallskip\noindent
{\bf Proposition:} There is a Quillen model structure on NSet with

\item{a)} weak equivalences ${\cal W}$ given by the Acyclic NSet maps,

\item{b)} fibrations ${\cal F}$ given by the Surjecting NSet maps,

\item{c)} cofibrations ${\cal C}={}^\dagger\underline{\cal F}$, where 
$\underline{\cal F}={\cal W}\cap{\cal F}$.

Moreover, the acyclic cofibrations $\underline{\cal C}$ are given by the Whiskering NSet maps.

\smallskip\noindent
{\bf Proof:} We create a model structure on NSet by applying 
the Transport Theorem to the adjunction
$$F :{\rm Gph}\rightleftharpoons\hbox{NSet}: G$$
We must check that $G ({\rm cell}(F J))\subseteq {\cal W}$. 
But $J$ contains just  the single graph morphism ${\bf s}:{\bf D}\to {\bf A}$.
We calculate that $F({\bf D})=(\pi_0({\bf P}\times {\bf D}),\sigma)=(N,\sigma)$ 
and $F({\bf A})=(\pi_0({\bf P}\times {\bf A}),\sigma)=(N,\sigma)$;
and the NSet map $F({\bf s}):(N,\sigma)\to(N,\sigma)$ is given by 
the successor function $\sigma:N\to N$.  Here are the details.  The
graph ${\bf P}\times {\bf D}$ has  no arcs and nodes $n$ for $n\geq 0$;
the graph ${\bf P}\times {\bf A}$ has nodes $(n,0)$ and $(n,1)$ for $n\geq 0$,
and an arc $(n,0)\to (n+1,1)$ for each  $n\geq 0$; and
${\bf s}:{\bf P}\times {\bf D}\to {\bf P}\times {\bf A}$ is given on nodes by 
${\bf s}(n)=(n,0)$. Then $\pi_0({\bf P}\times {\bf D})$ has elements $n$ for $n\geq 0$
and $\pi_0({\bf P}\times {\bf A})$ has elements $[n,1]$ for $n\geq 0$, but on components
we have ${\bf s}(n)=[n,0]=[n+1,1]$.
Let us show that $({\rm cell}(F J))$ is given by the Whiskering NSet maps.
Recall that the Whiskerings in Gph are attaching rooted trees, with arcs leaving the root.  
Such a  rooted tree is not in NGph, because the root is a node with no arcs entering.
Let us say that a {\it taprooted tree} is a rooted tree with a copy of the N-graph ${\bf N}$
(an infinite sequence of arcs and nodes leading into the node $0$)
attached to it by identifying its root with the $0$ node.
The NSet maps in ${\rm cell}(FJ)$ all come from attaching ``taprooted forests'',
and these are the Whiskering NSet maps.
Since $G$ preserves limits and colimits, $G({\rm cell}(F J))={\rm cell}(GF J)$. 
It follows that the  every graph morphism in $G({\rm cell}(F J))$ is a Whiskering, 
which is thus an Acyclic, and the hypothesis of the Transport Theorem is met.
Thus ${\cal F}$ and $\underline{\cal F}$ in NSet are defined in terms of the functor $G$,
which is the inclusion of NSet as the full subcategory NGph in Gph.
If $X=G(S,\sigma)$ then $C_n(S,\sigma)=C_n(X)$. 
Thus the weak equivalences in NSet are the Acyclic NSet maps.
Also, we have $X(x,*) = \sigma^{-1}(x)$,
so $G(f):G(S,\sigma)\to G(S',\sigma)$ is a Surjecting graph morphism
if and only if $f$ is a Surjecting NSet map.  
Thus the fibrations in NSet are the Surjecting NSet maps.
QED.

\bigskip
We define the Acyclic ZSet maps and the Surjecting ZSet maps by exactly copying the 
definitions used for NSet.
But these definitions simplify quite a bit in the category ZSet.

\smallskip\noindent
{\bf Proposition:} All ZSet maps are Surjecting.   

\smallskip\noindent
{\bf Proof:} In any ZSet $(S,\sigma)$, the function $\sigma$ is invertible, 
so $\sigma^{-1}(x)$ has exactly one element for every $x\in S$. Consider any ZSet
map $f:(S,\sigma)\to (T,\sigma)$. For every $x\in S$ we restrict $f$ to give 
$f:\sigma^{-1}(x)\to \sigma^{-1}(f(x))$, and any function between one element sets is surjective;
thus every ZSet map is Surjecting. QED. 

We can use the following definition to
describe the Acyclic ZSet maps. 
For any ZSet $(S,\sigma)$, let 
$\jj{(S,\sigma)}=\{x\in S:\exists n>0, \sigma^n(x)=x\}$. 

We may call $\jj(S,\sigma)$ the {\it periodic part} 
of the Z-set $(S,\sigma)$.
Any ZSet map  $f:(S,\sigma)\to (T,\sigma)$ restricts to give 
$\jj(f):\jj(S,\sigma)\to \jj(T,\sigma)$, since if
$\sigma^n(x)=x$ for some $x\in S$, then $\sigma^n(f(x))=f(x)$ in $T$.

\smallskip\noindent
{\bf Proposition:} A ZSet map $f$ is Acyclic if and only if $\jj(f)$ is a bijection. 

\smallskip\noindent
{\bf Proof:} 
Suppose that $f$ is Acyclic.  We want to show that $\jj(f)$ is a bijection. 
Certainly, 
$\jj(f)$ is a surjection, since for every $y\in \jj{T}$ we have 
$y\in C_n(T)$ for some $n>0$,
and $C_n(S)\to C_n(T)$ is bijective by assumption.  So there is a unique 
$x\in C_n(S)$ with $f(x)=y$.
Suppose that $\jj(f)$ is not an injection; then there exists some 
$y\in \jj{T}$ with more 
than one preimage in $\jj{S}$. We know that $y\in C_n(T)$ for some $n>0$, 
let $n$ be the smallest such.
Then there is a unique $x\in C_n(S)$ with $f(x)=y$.  
We have assumed there is another element $x'\in\jj{S}$ with
$f(x')=y$.  So $x'\notin C_n(S)$, and $x'\in C_m(S)$ for some $m>0$ with $m\neq n$.  
But then $f(x')=y$ must be in $C_m(T)$, and it follows that $m$ is a proper multiple of $n$.
Then we have $x,x'\in C_m(S)$, both mapping to $y \in C_m(T)$.  But $C_m(f)$ is a bijection. 
Contradiction. QED.

\smallskip\noindent
{\bf Corollary:} There is a Quillen model structure on ZSet with

\item{a)} the Acyclic ZSet maps as the weak equivalences ${\cal W}$,

\item{b)} all ZSet maps as the fibrations ${\cal F}$,

\item{c)} cofibrations ${\cal C}={}^\dagger{\underline{\cal F}}={}^\dagger{\cal W}$.

\smallskip\noindent
{\bf Proof:} Consider the adjoint functors
$$F :\hbox{NSet }\rightleftharpoons\hbox{ZSet}: G \quad{\rm or}\quad
F :\hbox{Gph}\rightleftharpoons \hbox{ZSet}: G$$
We use these to transport our model structure on Gph to  ZSet.
Note that $F({\bf s}):(Z,\sigma)\to (Z,\sigma)$
is the successor function $\sigma:Z\to Z$, which is an isomorphism.
This implies that every morphism in ${\rm cell}(F J)$ is an isomorphism,
so the hypothesis for the Transport Theorem is satisfied.  
Note that $\underline{\cal F}={\cal W}\cap{\cal F}={\cal W}$, since 
${\cal F}$ is all ZSet maps; thus ${\cal C}={}^\dagger{\cal W}$. QED.

We may summarize the above model structure on ZSet by
$$(\underline {\cal C},{\cal F})=({\rm iso},{\rm all})\quad{\rm and}\quad
({\cal C},\underline{\cal F})=({\cal C},{\cal W}).$$

\smallskip
Note that not every Acyclic ZSet map is an isomorphism.
For instance, any set $A$ gives a ZSet $A\times Z=\sum_{a\in A}Z$
(this is just viewing the set $A$ as a Z-set with trivial action, 
and taking the product of Z-sets).
If $A$ and $B$ are sets then any function $f:A\to B$ gives an ZSet map 
$f\times Z:A\times Z\to B\times Z$ by $(f\times Z)(a,n)=(f(a),n)$ for $a\in A$ and $n\in Z$.
Then $f\times Z$ is always Acyclic, but is usually not an isomorphism.

Let us say that an element $x$ in a Z-set $(S,\sigma)$ is {\it free} 
when $\sigma^n(x)=x$ implies $n=0$.  
The free part is not functorial, which seems dangerous.
However, we may say that a ZSet map $f:(S,\sigma)\to (T,\sigma)$ 
{\it maps the free elements bijectively} when $y\in T$ is free if and only 
if there exists a unique $x\in S$ with $f(x)=y$.  

\smallskip\noindent
{\bf Claim:} A ZSet map $f$ is a cofibration if and only if 
$f$ maps the free elements bijectively.

\beginsection ${\cal x}$3. Fibrant graphs and  cofibrant graphs.

A Quillen model structure on a category determines some important classes of objects there.
These are the fibrant, cofibrant, and fibrant-cofibrant objects.
In section 4 we will describe how they help to establish a well-behaved theory 
of homotopy classes of morphisms in the category.
In this section we investigate these notions for our model structure on Gph,
and describe a functor to ``replace'' any graph by a related cofibrant graph.

\smallskip\noindent 
{\bf Definition:}
Let $({\cal C},{\cal W},{\cal F})$ be
a model structure on a category ${\cal E}$, and let $X$ be an object in  ${\cal E}$.
We say that $X$ is {\it fibrant} when $X\to 1$ is in ${\cal F}$ (where $1$ is a terminal object);
we say that $X$ is {\it cofibrant} when $0\to X$ is in ${\cal C}$ (where $0$ is an initial object).
We say that $X$ is {\it fibrant-cofibrant} when it is both fibrant and cofibrant.

Let us see how this works in our model structures on Gph, NSet, and ZSet.
We start by introducing some terminology specialized to these different categories.

\item{*} For any graph $X$, a {\it dead-end} in $X$ is a node with no arc leaving it.

\item{*} For any N-set $(S,\sigma)$, 
the {\it trajectory} $N(x)$ of any element $x$ in $S$
is the set $\{\sigma^n(x): n\geq 0\}$.
We may say that an element $x$ is {\it periodic} when $\sigma^n(x)=x$ for some $n>0$. 
We  may say that $x$ is {\it eventually periodic} when $x$ has finite trajectory, since
$x$ has finite trajectory if and only if 
$\sigma^k(x)=\sigma^{n+k}(x)$ for some $n$ and $k$ (so that $\sigma^k(x)$ is periodic).

\item{*} For any Z-set $(S,\sigma)$,
the {\it orbit} $Z(x)$ for $x$ in $S$ is the set of elements $\sigma^n(x)$ 
as $n$ ranges over the integers. This definition makes sense 
for a Z-set since then the function $\sigma$ is invertible.

Note that we define the trajectory of an element in any N-set,
but the orbit of an element only makes sense in a Z-set, since the definition involves
the inverse function of $\sigma$.
An element in a Z-set is periodic iff it has a finite orbit.
In a Z-set, every element is either periodic (finite orbit) or free (infinite orbit).  
But in an N-set, an element with a finite trajectory may fail to be periodic.

\smallskip\noindent 
{\bf Proposition (Gph):}
A graph $X$ is fibrant if and only if $X$ has no dead-ends.
A graph $X$ is cofibrant if and only if $X$ is a disjoint union of whiskered finite-cycle graphs.
A graph $X$ is fibrant-cofibrant if and only if 
$X$ is a disjoint union of whiskered finite-cycle graphs 
in which the whiskers have no dead-ends.

\smallskip\noindent 
{\bf Proof:}  A graph $X$ is fibrant if and only if the morphism $X\to 1$ is Surjecting;
but this is true if and only if $X$ has at least one arc leaving each of its nodes.
A graph $X$ is cofibrant iff $0\to X$ is in ${\cal C}={}^\dagger(I^\dagger)$,
the class of retracts (in the morphisms category) of morphisms in ${\rm cell}(I)$.
If graph $X$ is a disjoint union of whiskered finite-cycle graphs, 
then $0\to X$ is in ${\rm cell}(I)$;
and every cofibrant graph is of this form, because
the only way to get the empty graph
$0$ as domain in ${\rm cell}(I)$ is to use a transfinite composition of
pushouts of the ${\bf i}_n$, and taking retracts can't introduce any
new morphisms with domain an empty graph.  
The description of fibrant-cofibrant graphs follows. QED.

\smallskip
It follows that the fibrant graphs are exactly those 
in which every path can be continued forever; this fits well with the terminology
``no dead-ends''. This can also be expressed by saying that any path 
${\bf P}_n\to X$ (for any length $n\geq 0$)
can be extended to an infinite path ${\bf P}\to X$.
Such graphs are convenient for the study of
symbolic dynamics and $C^*$ algebras (see Lind and Marcus [1995], for instance).

For example,  ${\bf Z}$ is not a cofibrant graph.
If $C$ is a cycle graph and $C\to CW$ adjoins some
whiskers with no dead-ends, then $C\to CW$ and $CW\to C$ are weak equivalences,
and both $C$ and $CW$ are fibrant-cofibrant.
Note that every cofibrant graph is in NGph, but not every fibrant graph is in NGph.

\smallskip
Recall that the category NSet is equivalent to the (reflective-coreflective) 
full subcategory NGph of Gph;
so we may carry-over some graph terminology in discussing N-sets.

\smallskip\noindent 
{\bf Proposition (NSet):}
An N-set $(S,\sigma)$ is fibrant if and only if $\sigma$ is surjective.
An N-set $(S,\sigma)$ is cofibrant if and only if each element in $S$ has finite trajectory.

\smallskip\noindent
{\bf Proof:}
The graph corresponding to an N-set $(S,\sigma)$ has a dead-end if and only if 
$\sigma$ is not surjective.
The graph corresponding to an N-set $(S,\sigma)$ is 
a disjoint union of whiskered finite-cycle graphs if and only if
each element in $S$ is eventually periodic, in that each trajectory in $S$ is finite. QED.

\smallskip\noindent 
{\bf Proposition (ZSet):} Every Z-set is fibrant.
A Z-set $(S,\sigma)$ is cofibrant if and only if each element in $S$ has finite orbit.

\smallskip\noindent
{\bf Proof:} Every ZSet map is a fibration.
The graph corresponding to a Z-set is cofibrant 
if and only if all of its connected components are finite-cycle graphs. QED.

\smallskip\noindent
{\bf Definition:} 
A {\it cofibrant replacement} for an object $X$ a model category ${\cal E}$ is a 
weak equivalence $f:X'\to X$ where $X'$ is cofibrant.
We will say that $f$ is a {\it full cofibrant replacement} if also $f\in \underline{\cal F}$.
Dually, a {\it fibrant replacement} of an object $Y$ is 
a weak equivalence $f:Y\to Y'$ where $Y'$ is fibrant; and
$f$ is a {\it full fibrant replacement} if also $f\in \underline{\cal C}$.

Each object $X$ in ${\cal E}$ has at least one
cofibrant replacement, since $0\to X$
has a $({\cal C},\underline{\cal F})$ factorization with
$0\to X'$ in ${\cal C}$ and $f:X'\to X$ in $\underline{\cal F}$; so this is 
actually a full cofibrant replacement.
Dually, we have the existence of (full) fibrant replacements.
If $X'$ is cofibrant and $g:X'\to X''$ is a full fibrant replacement, then
$X''$ is fibrant-cofibrant, since 
${\cal C}={}^\dagger\underline{\cal F}$ is closed under composition
($0\to X'$ in ${\cal C}$ and $g$ in $\underline{\cal C}$
implies $0\to X''$ in ${\cal C}$).
Thus any full fibrant replacement of a cofibrant object is fibrant-cofibrant,
Dually, a full cofibrant replacement of a fibrant object is fibrant-cofibrant.

\smallskip
Let us define a special cofibrant replacement for our model structure on Gph.
Recall, from the end of section 1, the adjoint pair
$$G:\hbox{ZSet}\rightleftharpoons \hbox{Gph}:H$$
Here $H$ is given by specifying the action on the set 
$H(X)=[{\bf Z},X]$ of morphisms from the line graph ${\bf Z}$ to the graph $X$; and
$G$ can be thought of as a Cayley graph construction, 
which is equivalent to the inclusion of ZGph as a subcategory of Gph.
The adjoint pair $(G,H)$ has counit $G(H(X))\to X$.
Recall that for any Z-set $(S,\sigma)$ we have defined
$\jj(S,\sigma)$ as the set of all elements $x\in S$ such that $\sigma^n(x)=x$
for some $n>0$.
Since $\sigma$ carries $\jj{S}$
into itself, we have a functor $\jj:\hbox{ZSet}\to\hbox{ZSet}$.
Then we may define $c(X)=G(\jj{H(X)})$.
We may refer to $c(X)$ as the {\it cycle resolution} of $X$.
For example, $c({\bf C}_n)={\bf C}_n$, and $c(X)=0$ if $X$ is an acyclic graph. 

Applying $G$ to $\jj{H(X)}\subseteq H(X)$ gives a natural graph morphism
$$c(X)=G(\jj{H(X)})\to G(H(X))\to X.$$
We have the following:

\smallskip\noindent
{\bf Proposition.}  For every graph $X$, the graph morphism $c(X)\to X$ is 
a cofibrant replacement.  It is not, in general, a full cofibrant replacement.

\smallskip\noindent
{\bf Proof:} We must show that $c(X)$ is cofibrant and that 
$c(X)\to X$ is an Acyclic graph morphism.
Note that $c(X)$ is always isomorphic to a disjoint union of finite cycle graphs,
and is thus a cofibrant graph.  Clearly $\jj{H(X)}\to H(X)$ is an
Acyclic ZSet map.  So applying $G$ to it gives an Acyclic graph morphism.
In other words, $c(X)$ is a natural subgraph (in fact, a summand) of $GH(X)$;
we merely omit all the ${\bf Z}$ components.
So it remains to show that $h:G(H(X))\to X$ is an Acyclic graph morphism;
in other words, that $C_*(h)$ is a bijection.
 Any ${\bf C}_n\to X$ is the image under $C_n(h)$ of the graph
morphism ${\bf C}_n=G(H({\bf C}_n)\to G(H(X))$, so $C_*(h)$ is surjective.
Conversely, any $\alpha:{\bf C}_n\to GH(X)$ gives $h\circ\alpha:{\bf C}_n\to X$.
Applying the functor $H$ gives $H(\alpha):H({\bf C}_n)\to HGH(X)$ 
and  $H(h\circ\alpha):H({\bf C}_n)\to H(X)$. But $HGH(X)=H(X)$,
since $HG$ is the identity on any Z-set.  Making this identification,
we have $H(\alpha)=H(h\circ\alpha)$.  Applying the functor $G$ gives
$GH(\alpha)=GH(h\circ\alpha)$, with $GH(\alpha):GH({\bf C}_n)\to GHGH(X)$.
But  $GH({\bf C}_n)={\bf C}_n$ and $GHGH(X)=GH(X)$ (as above);
after making these identifications, we have $GH(\alpha)=\alpha$.
It follows that $C_*(h)$ is injective, so $C_*(h)$ is a bijection.
To see that $c(X)\to X$ is not in general a full cofibrant replacement, we 
consider the graph $X$ with two nodes, $0$ and $1$, and two arcs,
from $\ell:0\to 0$ and $a:0\to 1$. Then $c(X)=1$ and the cofibrant replacement
$c(X)\to X$ is not Surjecting, and is thus not in $\underline{\cal F}$.
 QED

\smallskip\noindent
{\bf Definition:}
Let cZSet denote the full subcategory of ZSet whose objects are those 
with every element is periodic (those with no free elements).
Let $\ii$ denote the inclusion of 
cZSet as full subcategory of ZSet.
We may reinterpret $\jj$ (described above) as the left adjoint in the adjoint pair of functors
$$\ii :\hbox{cZSet}\rightleftharpoons \hbox{ZSet}: \jj.$$

The functor $\ii\circ\jj$ from ZSet to ZSet (with image cZSet)
is a ``comonad'' on ZSet, and cZset is isomorphic
to the topos of ``coactions'' for this comonad.
This exhibits cZSet as
a ``quotient topos'' of ZSet. See Mac Lane and Moerdijk [1994]  for a
discussion of these concepts.

Let ${G}\ii=G\circ\ii$ and $\jj{H}=\jj\circ H$, with adjoint pair
$${G}\ii :\hbox{cZSet}\rightleftharpoons\hbox{Gph}:\jj{H}$$
which results from the composition of the two adjoint pairs
$$\ii :\hbox{cZSet}\rightleftharpoons \hbox{ZSet}: \jj
\quad{\rm and}\quad G :\hbox{ZSet}\rightleftharpoons \hbox{Gph}: H.$$
Then we may interpret the cofibrant replacement functor $c:\hbox{Gph}\to\hbox{Gph}$ as
$c={G}\ii\circ \jj{H}$, which is the counit of the adjoint pair $(G\ii,\jj H)$.

We could also consider the adjoint functors 
$$\hbox{cZSet} \rightleftharpoons\hbox{ZSet}\rightleftharpoons
\hbox{NSet}\rightleftharpoons\hbox{Gph}.$$

\beginsection ${\cal x}$4. Homotopy categories.

Quillen [1967] introduced model categories
as a framework for defining and working with homotopy categories.
We discuss homotopy functors in general, then Quillen's definition
of the homotopy category as a category of fractions, and then
Quillen adjunctions and equivalences.
 
Suppose that we are given a model structure $({\cal C},{\cal W},{\cal F})$ 
on a category ${\cal E}$.
Recall that the morphisms in ${\cal W}$ are called weak equivalences. 
In the homotopy category, these should all become isomorphisms.  Let us sneak up on this idea.
We will say that a functor with domain ${\cal E}$
is {\it homotopy functor} when it takes every $f\in{\cal W}$ to an isomorphism.
We want to understand the homotopy functors for our model structure on Gph.
Consider functors from Gph to Set, for example.  Recall the adjoint triple $(F,G,H)$,
with $F(X)=\pi_0(X)$ and $H(X)=C_1(X)$.  Then

\smallskip\noindent{\bf Proposition:}
$H$ is a homotopy functor, and $F$ is not a homotopy functor.

\smallskip\noindent{\bf Proof:} 
The functor $H$ is clearly a homotopy functor, since $C_1(f)$ is a bijection 
for every graph morphism in ${\cal W}$
For the second part, the following example shows that $F$ is not a homotopy functor: 
the graph morphism $f:0\to 1$ is in ${\cal W}$ since $C_n(f)$ 
is the bijection $\emptyset\to \emptyset$ for all $n>0$.
But $\pi_0(f)$ is $\emptyset\to 1$, which is not a bijection of sets. QED

Recall the related adjoint triples $(F,G,H)$ relating Gph with NSet and with ZSet.
We will use subscripts to distinguish the cases.

\smallskip\noindent
{\bf Proposition:} The functors $F_N,H_N:{\rm Gph}\to {\rm NSet}$ 
and $F_Z,H_Z:{\rm Gph}\to {\rm ZSet}$ are not homotopy functors.

\smallskip\noindent
{\bf Proof:} Recall that the $F$ functors are reflecting.
The functor $\pi_0$ is a composition of the reflection functors
$$F:{\rm Gph}\to {\rm NSet}\to {\rm ZSet}\to {\rm Set}$$
It  follows that neither $F_{\rm N}$ nor $F_{\rm Z}$ is a homotopy functor, since if 
$F_1: {\rm Gph}\to {\cal A}$ is a homotopy functor,
and $F_2:{\cal A}\to {\cal B}$ is any functor, then
$F_2\circ F_1:{\rm Gph}\to {\cal B}$ must be a homotopy functor. 
The following example shows that neither of the $H$ functors is a homotopy functor: 
the graph morphism $f:0\to {\rm Z}$ is in ${\cal W}$ since $C_n(f)$ 
is the bijection $\emptyset\to \emptyset$ for all $n>0$.
But $H_{\rm N}(f)=H_{\rm Z}(f)$ is $\emptyset\to {\rm Z}$, which is not a bijection
(we used $H_{\rm N}(X)=[{\rm N},X]$ and 
$H_{\rm Z}(X)=[{\rm Z},X]$ to carry out this calculation). QED.

In the discussion below we will show that $jH:{\rm Gph}\to {\rm cZSet}$ is a homotopy functor.
This result underlies our calculation of the homotopical algebra of graphs.

Quillen [1967] showed how to use a model structure to
avoid set theoretic difficulties in the construction of 
a ``category of fractions'' which
universally inverts the morphisms in ${\cal W}$
so that they become isomorphisms.  

More precisely, Quillen used a model structure $({\cal C},{\cal W},{\cal F})$
to describe a  particular category ${\rm Ho}({\cal E})$,
together with a functor $\gamma:{\cal E}\to {\rm Ho}({\cal E})$ which is
{\it initial} for the homotopy functors (${\cal W}$-inverting functors) on ${\cal E}$.
This means that $\gamma$ is a homotopy functor and that
any homotopy functor $\Phi:{\cal E}\to{\cal D}$ factors 
uniquely through $\gamma$, in that $\Phi=\Phi'\circ \gamma$ for a unique functor 
$\Phi':{\rm Ho}({\cal E})\to {\cal D}$.

For example, if we use the trivial model structure $({\rm all},{\rm iso},{\rm all})$ on ${\cal E}$,
then ${\rm Ho}({\cal E})$ is isomorphic to ${\cal E}$.  This will apply to our model structure on
cZSet, and will help us describe ${\rm Ho}({\rm Gph})$.

In Quillen's description, the objects of the category ${\rm Ho}({\cal E})$ 
are the objects of ${\cal E}$. 
It follows that this universal definition determines ${\rm Ho}({\cal E})$ 
up to isomorphism of categories.
The category ${\rm Ho}({\cal E})$ is called the {\it homotopy category} for the model structure. 

The universal definition of ${\rm Ho}({\cal E})$ does not involve the fibrations and cofibrations, 
but these are used in Quillen's description of 
the set of morphisms from $X$ to $Y$ in ${\rm Ho}({\cal E})$, 
for objects $X$ and $Y$ in ${\cal E}$. 
We may denote this homotopy morphism set by ${\rm Ho}(X,Y)$.

Here is a sketch of Quillen's description of ${\rm Ho}(X,Y)$ 
for any objects $X$ and $Y$ in ${\cal E}$.
It uses the fibrations ${\cal F}$ and the cofibrations ${\cal C}$ as a kind of ``scaffolding''.
Suppose that $X'\to X$ and $Y'\to Y$ are full cofibrant replacements,
and that $X'\to X''$  and $Y'\to Y''$ are full fibrant replacements.  
(note that $X''$ and $Y''$ are objects which are both fibrant and cofibrant).
Then any $f:X\to Y$ can be factored by $f':X'\to Y'$, which can be
factored by some $f'':X''\to Y''$.
Quillen defines a ``homotopy''  equivalence relation $\ \sim\ $ on ${\cal E}(X'',Y'')$,
and uses the fibrant, cofibrant scaffolding to 
formally define ${\rm Ho}(X,Y)={\cal E}(X'',Y'')/\sim$.
Quillen shows that this definition supports a well-defined composition 
(which is independent of the choice of scaffolding),
and that this gives the category ${\rm Ho}({\cal E})$ with
functor $\gamma:{\cal E}\to {\rm Ho}({\cal E})$.

The functor $\gamma:{\cal E}\to {\rm Ho}({\cal E})$ gives a function
$\gamma:{\cal E}(X,Y)\to {\rm Ho}(X,Y)$;
we may denote $\gamma(f)$ by $[f]$.
However, general morphisms in ${\rm Ho}({\cal E})$ are
zig-zag compositions of homotopy classes of morphisms in ${\cal E}$;
the function $\gamma$ is not always surjective. 

We say that two objects in ${\cal E}$ are said to be {\it homotopy-equivalent} when
they become isomorphic in ${\rm Ho}({\cal E})$.

Suppose that ${\cal E}_1$ and ${\cal E}_2$ are Quillen model categories,
with model structure $({\cal C}_i,{\cal W}_i,{\cal F}_i)$ for ${\cal E}_i$ ($i=1,2$).
If $F:{\cal E}_1\to{\cal E}_2$
satisfies $F({\cal W}_1)\subseteq {\cal W}_2$, then 
$\gamma\circ F: {\cal E}_1\to{\rm Ho}({\cal E}_2)$ is a homotopy functor, and thus
factors through a unique ${\rm Ho}({\cal E}_1)\to {\rm Ho}({\cal E}_2)$.

Most functors we consider don't satisfy such a strong condition.
But Quillen developed a notion of {\it derived functor} 
suitable for homotopical algebra.
Suppose that $(L,R)$ is an adjoint pair of functors
$$L :{\cal E}_1\rightleftharpoons {\cal E}_2: R$$
between Quillen model categories ${\cal E}_1$ and ${\cal E}_2$. 
We say that $(L,R)$ is a {\it Quillen adjunction} when we have
$L({\cal C}_1)\subseteq {\cal C}_2$ and 
$L(\underline{\cal C}_1)\subseteq \underline{\cal C}_2$.
It turns out to be equivalent to have $R({\cal F}_2)\subseteq {\cal F}_1$ and 
$R(\underline{\cal F}_2)\subseteq \underline{\cal F}_1$ (see Hovey [1999], for instance).

A Quillen adjunction $L :{\cal E}_1\rightleftharpoons {\cal E}_2: R$ 
between model categories leads to an adjunction between
the respective homotopy categories, by means of derived functors. 
More precisely, Quillen described  a {\it (total) left derived functor} 
$L'$ associated to $L$, and a  {\it (total) right derived functor} $R'$ associated to $R$, 
giving an adjunction
$$L' :{\rm Ho}({\cal E}_1)\rightleftharpoons {\rm Ho}({\cal E}_2): R'$$

Here is a sketch of $L'$.
Suppose we choose a full cofibrant replacement $X'\to X$ for each object $X$. 
We can define $L'(X)=L(X')$;
and we can define $L'([f])=[f']$ for any $f:X\to Y$, where
$X'\to X$ and $Y'\to Y$ are full cofibrant replacements 
and $f':X'\to Y'$ is a lifting of $f$.  This definition of $L'$ 
extends uniquely to all morphisms in ${\rm Ho}({\cal E}_1)$.

The functor $L'$ comes with a natural transformation $\epsilon: L'\circ \gamma\to \gamma \circ L$,
and it is final (closest on the left) among all such factorizations through $\gamma$ 
(see Proposition 1 of chapter I.4 in Quillen [1967]).
This means that there is a unique $L''\to L'$ for any natural transformation 
$\epsilon': L''\circ \gamma\to \gamma\circ L$.
This universal condition determines the left derived functor $L'$ 
up to natural isomorphism of functors.

The description of $R'$ is dual to this, using full fibrant replacements, 
and has a dual universal property
(initial, or closest {\it on the right} among all factorizations through $\gamma$).  
See chapter I.4 in Quillen [1967] for more details.

Here are some examples of derived functors on Gph.
Recall the adjoint pair $F :{\rm Gph}\rightleftharpoons {\rm Set}: G$
with $F(X)=\pi_0(X)$ and $G(S)=\sum_S 1$.  Consider our model structure on Gph
and the trivial model structure on Set.

\smallskip\noindent{\bf Proposition:} $(F,G)$ is a Quillen adjunction,
with derived adjunction $F' :{\rm Ho(Gph)}\rightleftharpoons {\rm Ho(Set)}: G'$
having $G'(S)=G(S)$ and $F'(X)=\pi_0(c(X))$.

\smallskip\noindent{\bf Proof:}  Recall the 
trivial model structure $({\rm all},{\rm iso},{\rm all})$ on Set,
and let $({\cal C},{\cal W},{\cal F})$ be our model structure on Gph.
It is easy to see that $G(f)\in{\cal C}$ for any function $f$, 
and $G(h)\in{\cal C}$ for any bijection $h$;
so $(F,G)$ satisfy the Quillen conditions.  In fact, $G$ is clearly a homotopy functor, 
so we may take $G'\circ\gamma=G$.
We claim that we can take $F'(X)=F(c(X))$ in the left derived functor, despite the fact that
$c(X)\to X$ is not a {\it full} cofibrant replacement.  This will become easy to verify when
we have finished our calculation of Ho(Gph) by the end of this section. QED

It is not hard to show that the adjoint pairs
$$F :{\rm Gph}\rightleftharpoons {\rm NSet}: G \quad{\rm and}\quad
F :{\rm Gph}\rightleftharpoons {\rm ZSet}: G$$
are also Quillen adjunctions; but we don't need them here.
We find it more convenient to deal instead with adjunctions
$$G :{\rm ZSet}\rightleftharpoons {\rm Gph}: H \quad{\rm and}\quad
G \ii :{\rm ZSet}\rightleftharpoons {\rm Gph}: \jj H$$ 
Will we show that these are Quillen adjunctions and that the derived adjunctions are actually
equivalences of categories.  This will show that ${\rm Ho(Gph)}$ and ${\rm Ho(ZSet)}$ 
are both equivalent to the category cZSet.  Let us begin.

Recall that an equivalence of categories is just a special kind of adjunction.
A Quillen adjunction $(L,R)$ for which $(L',R')$ is an equivalence 
is called a {\it Quillen equivalence}.  We will use the following characterization
of Quillen equivalence.

\smallskip\noindent
{\bf Theorem:} A Quillen adjunction $(L,R)$ is a Quillen equivalence if and only if
for all cofibrant $X$ in ${\cal E}_1$ and all fibrant $Y$ in ${\cal E}_2$ we have
$LX\to Y$ in ${\cal W}_2$ if and only if $X\to RY$ in ${\cal W}_1$.

We start by applying  this to the following situation.  Recall the definition of cZSet
as the full subcategory of Z-sets in which every element is periodic,
as discussed in the previous section).
Consider the functors $\ii:{\rm cZSet}\to {\rm ZSet}$ and $\jj:{\rm ZSet}\to {\rm cZSet}$
where $i:{\rm cZSet}\to{\rm ZSet}$ is the inclusion of the subcategory, and
$\jj(S,\sigma)=(\jj{S},\sigma)$ (as discussed in the previous section).
Consider the trivial model structure on cZSet.  Recall the model structure on ZSet 
which we described
in section 2. 

\smallskip\noindent
{\bf Proposition:} The functors $\ii$ and $\jj$ are adjoint, and 
the adjoint pair $\ii :{\rm cZSet}\rightleftharpoons  {\rm ZSet}: \jj$ is a Quillen equivalence.
So  ${\rm Ho(ZSet)}$ and ${\rm Ho(cZSet)}$ are equivalent
by $\ii'$ and $\jj'$, and ${\rm Ho(cZSet)}$ is equivalent to cZSet.

\smallskip\noindent
{\bf Proof:}  First show the adjunction $$\ii :{\rm cZSet}\rightleftharpoons {\rm ZSet}: \jj$$
Then show that $\jj$ behaves correctly on fibrations and acyclic fibrations.
Then show that the Quillen equivalence condition is satisfied. QED

\smallskip\noindent
{\bf Proposition:} $(G \ii,\jj H)$ is a Quillen equivalence.
So the homotopy category ${\rm Ho(Gph)}$ is equivalent to the category cZSet.

\smallskip\noindent
{\bf Proof:} It is easy to check that $G\ii$ satisfies the Quillen adjunction condition 
(on fibrations and acyclic fibrations).
Then use $c(x)=G\ii\jj H(X)$ and $c(X)\to X$ is an Acyclic, for any graph $X$.
QED.

\beginsection ${\cal x}$5. Isospectral graphs.

It seems that our homotopy category of graphs fits well with algebraic graph theory
and other parts of combinatorics.  
Let us illustrate this by connecting the treatment of zeta series in Bisson and Tsemo [2008]
with that in Dress and Siebeneicher [1988] and [1989].  
 
There they work with Burnside rings of Z-sets and actions of profinite groups, and show
how this algebra is mirrored in theories of zeta series and Witt vectors.
 Recall that a Z-set is a set $S$ together with an invertible
 function $\sigma:S\to S$.  For example, the integers modulo $n$ form a Z-set by taking
  $\sigma(i)=i+1\ {\rm mod}\ n$; let us denote this Z-set by $Z/n$.   For Z-sets $S$ and $T$, let $[S,T]$ denote  the set of Z-set maps from $S$ to $T$. Let $Z(x)$ denote the orbit of an element $x$ in a Z-set.
 Recall that we say that an element $x$ in a Z-set is periodic
 when $Z(x)$ is finite.

 \smallskip\noindent
{\bf Definition (Dress and Siebeneicher):}  A Z-set $S$ is {\it essentially-finite} 
when $[Z/n,S]$ is finite for all $n>0$. The  {\it zeta series} of an essentially-finite Z-set $S$
 is defined by
$$Z_S(u)={\rm exp}(\sum_{n=1}^\infty c_n{u^n\over n}),$$
where $c_n$ is the cardinality of $[Z/n,S]$, for all $n>0$.
An {\it almost-finite} Z-set is an essentially-finite Z-set for which 
every element is periodic.

Recall that $\jj:{\rm ZSet}\to {\rm cZSet}$  is given by taking the periodic part of a Z-set,
so that $\jj(S)$ is the set of periodic elements in $S$. Note that if $S$ is essentially-finite then
$\jj(S)$ is almost-finite; in fact, $S$ and $\jj(S)$ have the same zeta series.
This follows from the fact that, by
definition, two essentially-finite Z-sets $S$ and $T$ have the same zeta series if and only if
$[Z/n,S]$ and $[Z/n,T]$ have the same cardinality for all $n>0$.  
But  the following result is noted in Dress and Siebeneicher [1989] (without proof).

\smallskip\noindent
{\bf Proposition:}  
Two  almost-finite Z-sets $S$ and $T$ have the same zeta series if and only if 
$S$ and $T$ are isomorphic as Z-sets.

\smallskip\noindent
{\bf Proof:}  Let $Z/(n)$ denote the finite cyclic group of order $n$.
If $S$ is an almost finite Z-set, then let $S_n=\{x: |Z(x)|=n\}$.  So $S$ is the disjoint union of the $S_n$,
and each $S_n$ is a finite and free $Z/(n)$-set.  Let  $s_n Z/n$ denote the 
sum ( disjoint union) of $s_n$ copies of  the Z-set $Z/n$, where $s_n=|S_n|/n$.
Thus $S_n$ is isomorphic to $s_n Z/n$ as Z-sets.
This shows that $S$ is isomorphic to $\sum_{n>0} s_n Z/n$.  It suffices to show that
the numbers $(c_n:n>0)$ determine the numbers $(s_n:n>0)$.  This is true because the
``triangular'' system of equations $c_n=\sum_{k|n} k s_k$ has a unique solution.  QED.

This result shows that assigning a zeta series to each almost-finite Z-set gives an
isomorphism between the Burnside ring of almost-finite Z-sets and the
universal Witt ring (with integer coefficients).  See Dress and Siebeneicher [1988] for details.

Let us lift some of these definitions to the category Gph, by using our functor
$H:{\rm Gph}\to {\rm ZSet}$; recall that $H(X)$ is the Z-set $[{\rm Z},X]$.

\smallskip\noindent
{\bf Definition:}
A graph $X$ is {\it essentially-finite} if $H(X)$ is an essentially-finite Z-set. 
The {\it zeta series} of an essentially-finite graph $X$
is the formal power series $Z_X(u)=Z_{H(X)}(u)$.
A graph $X$ is {\it almost-finite} when $H(X)$ is an almost-finite Z-set.

For example, the graph ${\bf Z}$ is  an essentially-finite graph which is not almost-finite.
Note that a graph $X$ is essentially finite if and only if $C_n(X)$ is finite for all $n>0$.
Also, a graph $X$ is essentially-finite if and only if $\jj H(X)$
is an almost-finite Z-set.

Let us review a few concepts from algebraic graph theory, with
terminology as in Bisson and Tsemo [2008].
A {\it finite graph} $X$ is one with finitely many nodes and arcs.
The {\it characteristic polynomial} of a finite graph $X$ is defined as 
$a(x)=\det(xI-A)$, the characteristic polynomial of the adjacency operator $A$ for $X$.  
If $X$ has $n$ nodes, then $a(x)$ is a monic polynomial of degree $n$, and
the {\it reversed characteristic polynomial} of $X$ is defined to be
$u^n a(u^{-1})=\det(I-uA)$.
The roots of the characteristic polynomial of $X$ form the {\it spectrum} of $X$
(the eigenvalues of the adjacency operator for $X$).  This motivates
the following.

\smallskip\noindent
{\bf Definition:}
Two finite graphs $X$ and $Y$ are {\it isospectral} 
if they have the same characteristic polynomial.  
Two finite graphs $X$ and $Y$ are  {\it almost-isospectral} if they have the same
reversed characteristic polynomial.

Loosely speaking, $X$ and $Y$ are almost-isospectral if and only if
they have the same non-zero eigenvalues,
since $u=z$ is a root of $\det(I-uA)$ if and only if $z\neq 0$ and $x=z^{-1}$ is a root of $\det(xI-A)$.

We give a proof in Bisson and Tsemo [2008] of the following (folk) result on directed graphs.

\smallskip\noindent
{\bf Proposition:} Two finite graphs $X$ and $Y$ are almost-isospectral 
if and only if  $Z_X(u)=Z_Y(u)$.

Note that every finite graph is almost-finite.  
We say that two graphs are {\it homotopy equivalent} when they become isomorphic in Ho(Gph).

\smallskip\noindent
{\bf Theorem:}  Two finite graphs $X$ and $Y$ are almost-isospectral if and only if 
they are homotopy-equivalent.

\smallskip\noindent
{\bf Proof:} If $X$ is a finite graph, then $Z_X(u)=Z_{H(X)}(u)=Z_{jH(X)}(u)$.
In section 4 we showed that $jH:\hbox{Gph}\to \hbox{cZSet}$ is a
homotopy functor.  So there is a unique functor $H'':\hbox{Ho(Gph)}\to \hbox{cZSet}$
such that $H''\circ \gamma_1=jH$ where $\gamma_1:\hbox{Gph}\to \hbox{Ho(Gph)}$.
In fact, $\gamma\circ H''$ is the total left derived functor of the Quillen equivalence $jH$,
where $\gamma_2:\hbox{cZSet}\to \hbox{Ho(cZSet)}$ is an isomorphism of categories.
So $H''$ is an equivalence of categories.

Suppose that $X$ and $Y$ are almost-isospectral.  Then $Z_X(u)=Z_Y(u)$.
So $Z_{jH(X)}(u)=Z_{jH(Y)}(u)$.  But $jH(X)$ and $jH(Y)$ are almost-finite Z-sets.
So by the Dress-Siebeneicher Proposition above, $jH(X)$ and $jH(Y)$ are isomorphic
in cZSet.  So $H''(\gamma_1(X))$ and $H''(\gamma_1(Y))$ are isomorphic
in cZSet.  But $H''$ is an equivalence of categories,  so $\gamma_1(X)$ and $\gamma_1(Y)$
are isomorphic in Ho(Gph); so $X$ and $Y$ are homotopy equivalent.

Conversely, suppose that $X$ and $Y$ are homotopy equivalent.  This means that $\gamma_1(X)$
and $\gamma_1(Y)$ are isomorphic.  So $H''(\gamma_1(X))=jH(X)$ is isomorphic to
$H''(\gamma_1(Y))=jH(Y)$ in cZSet.  So $jH(X)$ and $jH(Y)$ have the same zeta series.
So $X$ and $Y$ have the same zeta series, so $X$ and $Y$ are almost-isospectral. QED.

\smallskip\noindent
{\bf Corollary:} Two finite graphs have the same zeta series 
if and only if they are homotopy-equivalent in our model stucture for Gph.

\smallskip\noindent
{\bf Example:}  Consider the graph with vertices $0,1,2,3,4$ and arcs
$(0,i)$ and $(i,0)$ for $i=1,2,3,4$; we will call it the Cross.  
Let $U{\bf C}_4$ be the undirected cycle, with nodes the integers mod $4$,
with arcs $(i,i+1)$ and $(i,i-1)$ for all $i$ mod $4$,
and with source and target given by $s(i,j)=i$ and $t(i,j)=j$.
The characteristic polynomial of $U{\bf C}_4$ is $x^4-4x^2$;
the characteristic polynomial of the Cross is $x^5-4x^3$.
So they have the same reversed characteristic polynomial $1-4u^2$,
and thus the same  zeta series 
$$Z(u)=(1-4u^2)^{-1}=\sum_{n\geq 0} 2^{2n}u^{2n}={\rm exp}(\sum_{n>0} 2^{2n+1}{u^{2n}\over {2n}})$$
So $U{\bf C}_4$ and the Cross are almost-isospectral, 
and thus must be homotopically equivalent for our model structure for Gph.

We note in passing that the formula for $Z(u)$
says that there are no graph morphisms from an odd cycle to either graph,
and that there are exactly  $2^{2n+1}$ graph morphisms from ${\bf C}_{2n}$ to each graph
(since $c_{2n}=2^{2n+1}$).


\bigskip
\noindent{\bf References:}

\smallskip
\item{[2003]} C. Berger  and I. Moerdijk, Axiomatic homotopy theory for operads. Comment. Math. Helv. 78 (2003), no. 4, 805--831.

\smallskip
\item{[2008]} T. Bisson and A. Tsemo, A homotopical algebras of graphs
related to zeta series, Homology, Homotopy and its Applications, 10 (2008),
1-13. 

\smallskip
\item{[1995]} S.E. Crans, Quillen closed model structures for sheaves.  J.
Pure Appl. Algebra  101  (1995),  no. 1, 35--57.

\smallskip
\item{[1988]} A.W.M. Dress and C. Siebeneicher, The Burnside ring of
profinite groups and the Witt vector construction.  Adv. in Math.  70  (1988), 
no. 1, 87--132. 

\smallskip
\item{[1989]} A.W.M. Dress and C. Siebeneicher, The Burnside ring of the infinite
cyclic group and its relations to the necklace algebra, $\lambda$-rings, and the
universal ring of witt vectors. Adv. in Math., 78 (1989),1-41. 

\smallskip
\item{[2000]} E.J. Dubuc and C.S. de la Vega,  On the Galois theory of Grothendieck. in Colloquium on Homology and Representation Theory (Vaquer'as, 1998). Bol. Acad. Nac. Cienc. (C—rdoba) 65 (2000), 111--136. Also on ArXiv.

\smallskip
\item{[1995]} W.G. Dwyer and J. Spalinski, Homotopy theories and model categories,
in {\bf Handbook of Algebraic Topology}, Ed. I. M. James, 73-126, Elsevier,
1995. 

\smallskip
\item{[1999]} L. Fajstrup and J. Rosick\'y, A convenient category for directed
homotopy,Theory and Applications of Categories, 21 (2008), No. 1, 7-20. 

\smallskip
\item{[1972]} A. Grothendieck et al, Th\'eorie des topos et cohomologie \'etale des sch\'emas. 
Tome 1: Th\'eorie des topos, S\'eminaire de G\'eomŽtrie Alg\'ebrique du Bois-Marie 1963--1964 
(SGA 4). Dirig\'e par M. Artin, A. Grothendieck, et J. L. Verdier. 
Avec la collaboration de N. Bourbaki, P. Deligne et B. Saint-Donat. 
Lecture Notes in Mathematics, Vol. 269. Springer-Verlag, Berlin-New York, 1972.

\smallskip
\item{[1999]} M. Hovey, {\bf Model Categories}, Amer. Math. Soc., Providence, 1999.

\smallskip
\item{[2007]} A. Joyal and M. Tierney, Quasi-categories vs Segal spaces, 277-326 in
{\bf Categories in algebra, geometry and mathematical physics}, Contemp. Math.
431, Amer. Math. Soc., Providence, 2007. 

\smallskip
\item{[1989]} F.W. Lawvere, Qualitative distinctions between some toposes of
generalized graphs, 261-299 in {\bf Categories in computer science and logic
(Boulder 1987)}, Contemp. Math. 92, Amer. Math. Soc., Providence, 1989.

\smallskip
\item{[1997]} F.W. Lawvere and S.H. Schanuel, {\bf  Conceptual Mathematics: a first
introduction to categories}, Cambridge University Press, Cambridge, 1997. 

\smallskip
\item{[1995]} D. Lind and B. Marcus,  An introduction to symbolic dynamics and
coding. Cambridge University Press, Cambridge, 1995. 

\smallskip
\item{[1971]} S. Mac Lane,  Categories for the working mathematician. 
Graduate Texts in Mathematics, Vol. 5. Springer-Verlag, New York-Berlin, 1971 (Second
edition, 1998).

\smallskip
\item{[1994]} S. Mac Lane and I. Moerdijk, {\bf Sheaves in Geometry and Logic: a
first introduction to topos theory}, Universitext, Springer-Verlag, New York
(1994).

\smallskip
\item{[1967]} D.G. Quillen, Homotopical Algebra, Lecture Notes in Mathematics No.
43, Springer-Verlag, Berlin, 1967.

\end

\bigskip\bigskip\noindent
Terrence  Bisson, \quad bisson@canisius.edu

Department of Mathematics and Statistics, Canisius College, 2001 Main Street, Buffalo, NY 14216 USA

\bigskip\noindent
Aristide Tsemo, \quad tsemo58@yahoo.ca

College Boreal, 351 Carlaw Avenue, Toronto, Ontario M4K 3M2 Canada

MR codes:

05C20: directed graphs
05C20: cycles
18F20: presheaves
55U35: abstract and axiomatic homotopy theory

For our first graph paper we used these:
\classification{05C20,18G55, 55U35.}

\end